\documentclass[hidelinks,conference]{IEEEtran}
\usepackage{ifpdf}
\usepackage{amsfonts}
\usepackage{latexsym}
\usepackage{graphicx}
\usepackage{epsfig}
\usepackage{amsbsy}
\usepackage{amsmath}
\usepackage{amssymb}
\usepackage{psfrag}
\usepackage{subfigure} 
\usepackage{epstopdf}
\usepackage{microtype}
\usepackage{hyperref}
\usepackage{cite}
\usepackage{color}
\usepackage[T1]{fontenc}
\usepackage{tabularx,ragged2e,booktabs,caption}
\newcolumntype{C}[1]{>{\Centering}m{#1}}

\usepackage{lmodern}

\usepackage{array}

\newcommand{\E}{{\mathbb E}}

\begin{document}

\title{Revenue maximization in an optical router node -- \\
allocation of service windows}

\author{\IEEEauthorblockN{Murtuza Ali Abidini\IEEEauthorrefmark{1}\IEEEauthorrefmark{2},
Onno Boxma\IEEEauthorrefmark{1},
Ton Koonen\IEEEauthorrefmark{2},and
Jacques Resing\IEEEauthorrefmark{1}}
\IEEEauthorblockA{\IEEEauthorrefmark{1}EURANDOM and Department of Mathematics and Computer Science}
\IEEEauthorblockA{\IEEEauthorrefmark{2}COBRA Institute and Department of Electrical Engineering}
Eindhoven University of Technology\\
P.O. Box 513, 5600MB Eindhoven, The Netherlands\\

\href{mailto:m.a.abidini@tue.nl}{m.a.abidini@tue.nl}, \href{mailto:o.j.boxma@tue.nl}{o.j.boxma@tue.nl},
\href{mailto:a.m.j.koonen@tue.nl}{a.m.j.koonen@tue.nl}, \href{mailto:j.a.c.resing@tue.nl}{j.a.c.resing@tue.nl}}

\maketitle

\begin{abstract}
In this paper we study a revenue maximization problem for optical routing nodes.
We model the routing node as a single server polling model with the aim to assign
visit periods (service windows) to the different stations (ports) such that
the mean profit per cycle is maximized. Under reasonable assumptions regarding
retrial and dropping probabilities of packets the optimization problem becomes
a separable concave resource allocation problem, which can be solved using existing algorithms.
\end{abstract}


\IEEEpeerreviewmaketitle

\section{Introduction}
The traffic routing in telecommunication networks has undergone a dramatic shift in the last decades due
to the changing nature of telecommunication services: from slowly-changing circuit-switched traffic routes
for traditional voice telephony to highly dynamic packet-switched traffic routes for internet traffic.
Hence also the demands on the nodes in the network have become much higher, not only regarding sheer
traffic volume but also regarding reconfiguration times. In fast bidirectional interactive communication
(such as machine-to-machine), it is important to minimize latency, as any delay occurring in the network
will reduce its throughput and deteriorate the Quality-of-Service experienced at the user, in particular
for smaller-sized packet communication. E.g., in a TCP/IP based link the throughput is approximately
inversely proportional to the round-trip time in the link, and proportional to the TCP window size
(see e.g. van Mieghem  \cite{Mieghem}, Ch. 5). 

\begin{figure}[h]
\hspace*{0.1in}
 \includegraphics[width=0.5\textwidth, height=0.35\textheight]{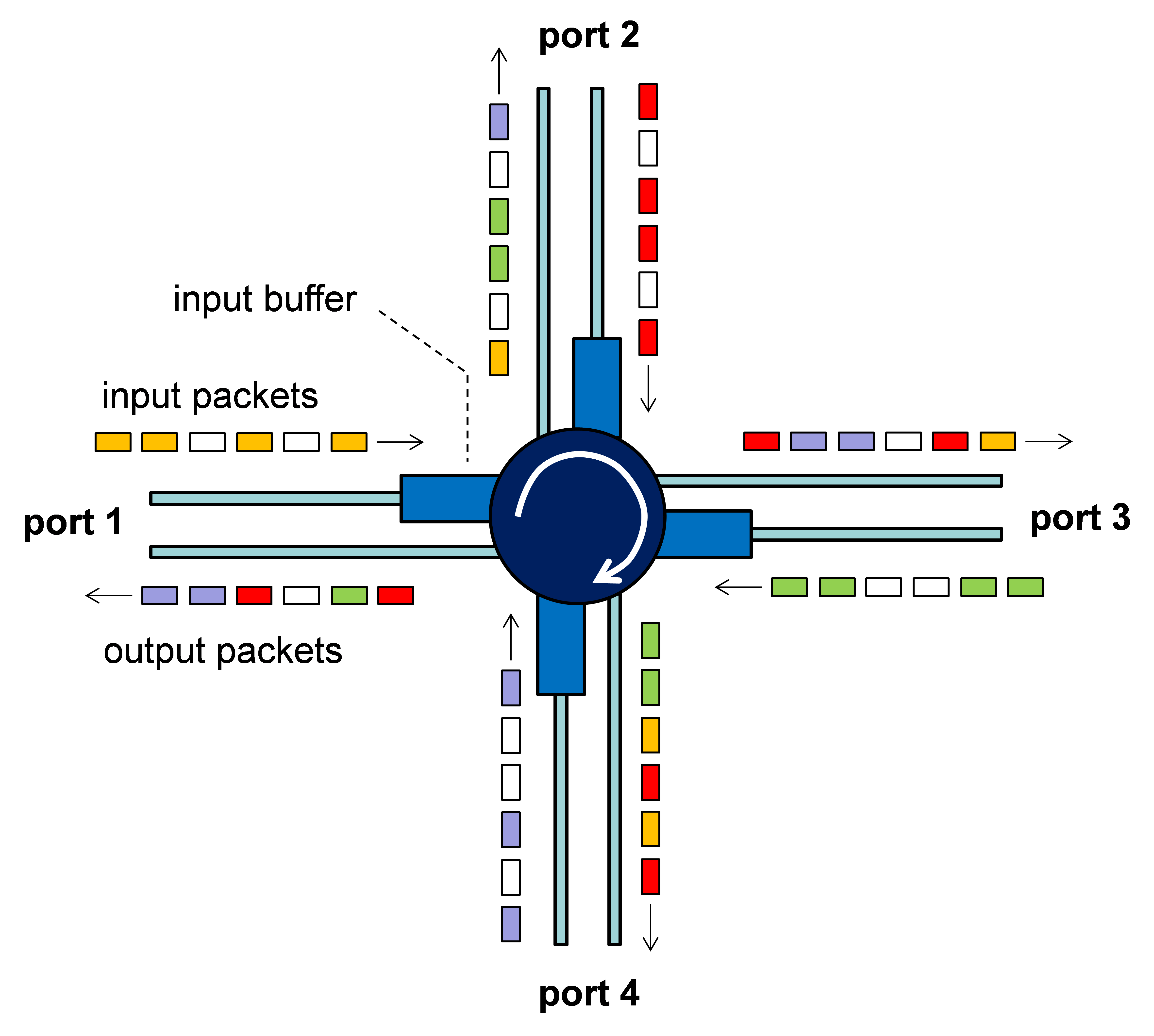} 
 \caption{Two-dimensional routing node (routing in time and in space)}
 \label{fig_12}
\end{figure}

In a telecommunication network, packets have to be routed from source to destination, passing through a sequence
of links and nodes. Packets from different sources are time-multiplexed and thus flow sequentially through the
network's links. When arriving at a routing node, they need to be queued in a buffer, where they need to wait
before they can be forwarded to the appropriate outgoing port of the node and travel further through the network;
see Fig. \ref{fig_12}. This store-and-forward procedure can cause a serious increase of the latency, and increasingly
so when the traffic load in the network grows. Hence the buffering processes need to be designed as efficiently as possible. 

\begin{figure}
 \includegraphics[width=0.5\textwidth, height=0.2\textheight]{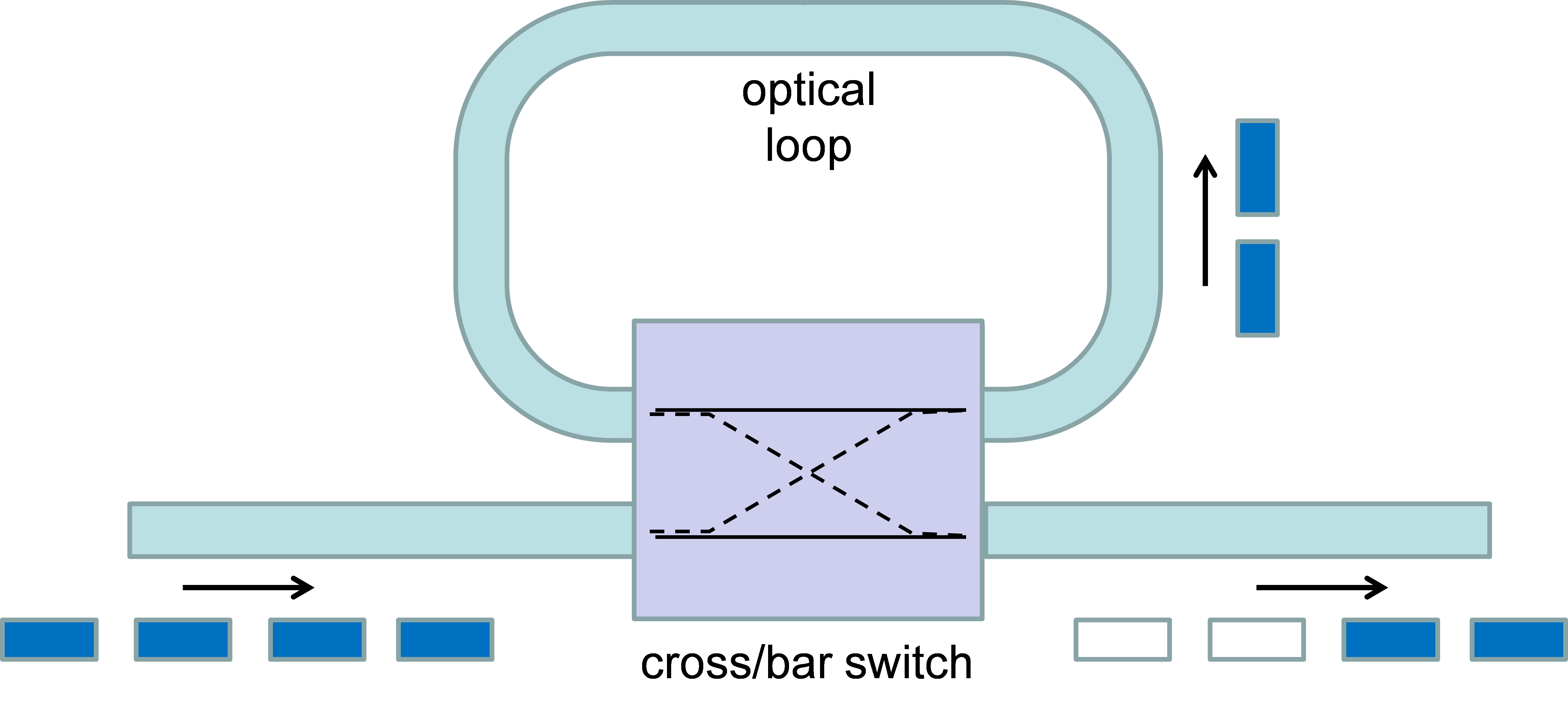} 
 \caption{Buffering packets in an optical delay loop}
 \label{fig_13}
\end{figure}

The performance analysis of optically-routed networks brings additional challenges with respect to the analysis
of networks which deploy electronic routing (see e.g. Maier \cite{Maier} and Rogiest \cite{Rogiest}). One of those
challenges is buffering, as in optical
networks it is difficult to store photons. Buffering in these networks is typically realised by
sending optical packets into fiber delay loops, i.e., letting them recirculate in a fiber loop and extracting
them after a certain number of circulations, as shown in Fig. \ref{fig_13}. Packets can be inserted into and
extracted from the delay loop by means of e.g. a cross/bar switch. This optical storage concept can be
modelled by so-called retrial queues.

In this paper, we report on the modelling of an optical routing node as a queueing system,
in which we aim to maximize its performance by a `revenue maximization' approach.
We develop a strategy to optimize the server times for the respective ports of the
node while taking into account the various retrial (buffering) times provided by the optical
delay loop buffering concept.

{\it Optical routing node represented by a polling model:}
We use a so-called polling model to study the performance of this routing node.
A polling model is a queueing model in which the server cyclically visits a number ($N$) of queues/stations,
serving customers at station~$i$ for a while and then switching to station~$i+1 ~{\rm mod} ~ N$, $i=1,2,\dots,N$.
Polling systems have been extensively studied in the literature. For example, various
different service disciplines (rules which describe the server's behaviour
while visiting a queue) have been considered, and models with and also without switchover times.
We refer to \cite{Takagi1,Takagi2} and \cite{Vishnevskii} 
for literature reviews and to \cite{Boon,Levy} and \cite{Takagi3} for overviews of the 
applicability of polling systems, with some emphasis on applications in communication systems.

{\it Choosing the polling times optimally:}
As a communication system typically works in frame time, which is fixed, we demand that the time it takes the server to complete one cycle of the
$N$~stations in the polling model is a given constant, $C$.
We want to assign fixed amounts of time $V_1,\dots,V_N$ to the visit periods (also called service windows) of stations $1,\dots,N$, such that
$\sum_{i=1}^N V_i = C - \sum_{i=1}^N S_i$, where $S_i$ denotes the time to switch to station~$i$, $i=1,2,\dots,N$.
If, say, $V_i$ is relatively small, then there is a relatively high probability that a packet in a retrial loop for station~$i$
does not retry during the visit period. Such packets may have to be dropped.
We assume that each served packet generates a profit, whereas each dropped packet incurs a loss to the system.
We allow at each station multiple customer types with different revenue/cost characteristics.

Our goal is to maximize the system revenue, under the above constraint regarding $\sum_{i=1}^N V_i$.
Under reasonable assumptions on the probability $p_i(V_i)$ that a packet in a retrial loop of station~$i$
actually retries during the visit period $V_i$,
and on the probability $q_i(V_i)$ that a packet is dropped when it fails to retry during $V_i$,
the revenue optimization problem turns out to be a so-called separable concave optimization problem.
This is a well-studied type of optimization problem, allowing for an efficient and quite insightful algorithm 
that yields the optimal solution.
We shall demonstrate the algorithm for some small examples.

{\it Organization of the paper:}
The model under consideration is described in
Section~\ref{sec:model}. In Section~\ref{sec:perform} we derive expressions for the mean numbers of customers in the stations/retrial loops at various epochs, and use these
to determine an expression for the mean revenue at each station per cycle.
In Section~\ref{sec:resource} we formulate the revenue optimization problem, show that
it indeed becomes a separable concave optimization problem under reasonable assumptions, and indicate
how it can be solved.
Section~\ref{sec:numeric} presents three numerical examples.
Section~\ref{sec:conclusion} contains conclusions and some suggestions for further research.

\section{Optical routing node model}
\label{sec:model}
Consider an optical routing node with $N$~ports (stations) to route packets and retrial loops to store packets. We represent it by a single
server polling model, i.e., a queueing model with a single server
which cyclically visits $N$~queues. Packets (also called customers in queueing terminology) of type~$j$, $j=1,\cdots,M$,
arrive at station~$i$ according to independent Poisson processes with rate $\lambda_{ij}.$ If at the time of arrival the station is being served
then the packet is instantaneously transmitted; else it enters a retrial loop.  In optical nodes the retrial
time is the delay produced by the fiber delay loop. We assume the retrial time to be random, because delay loops of various lengths may be used.  If, at the time of retrial, the
station is not in service then the packet again goes into a retrial loop and this process continues.

The server visits each station~$i$ for a fixed period of time $V_i$. During this period there may be two types of arrivals: (i) newly arriving packets, and (ii) 
packets which were in a retrial 
loop; we assume the latter retry during $V_i$ with probability $p_i(V_i)$. The server serves all these packets (new arrivals + retrials) instantaneously,
i.e., the service rate is assumed to be infinite.
At the end of the visit of station~$i$ each packet which still resides in a retrial loop of $i$ is dropped
with the same probability $q_i(V_i)$. Then the server moves to station~$i+1$ mod $N$ with a deterministic switchover time $S_{i+1 ~ {\rm mod} ~N}$. Hence the probability
that a packet in a retrial loop of station~$i$ leaves the system, either served during a visit at station~$i$ or dropped after a visit of station~$i$, is $r_i(V_i) := p_i(V_i)+q_i(V_i)-p_i(V_i)q_i(V_i)$.
Summarizing, 
\begin{itemize}
\item The customers of type~$j$ arrive at station~$i$ according to independent Poisson processes with rate $\lambda_{ij}$, $i= 1,
\cdots,N$ and $j= 1, \cdots,M$.
\item The length of a visit period at each station~$i$ is $V_i$.
\item The switchover time to station~$i$ is $S_i$.
\item The packets at station~$i$ retry during the visit period with probability $p_i(V_i)$.
\item After a visit at station~$i$, the packets in their retrial loops
are individually dropped with probability $q_{i}(V_i)$.
\item After a visit at station~$i$, each packet which resided in a retrial loop at the start of the visit
has left the system with probability $r_{i}(V_i)$.
\end{itemize}
\vspace{0.2in}
The motivation behind the model is as follows:
\begin{itemize}
 \item Since an optical routing node has multiple input ports we assume $N$ ports.
 \item Since the buffers used to store an optical packet are fiber delay loops we assume retrial loops.
 \item We consider a single-wavelength system in which only one port can transmit at a time,
 hence  we assume a single server with cyclic service.
 \item Since there can be more than one type of data at each port we assume that there are $M$ types of arrivals which are independent of each other.
 \item Since the server needs a positive amount of time to change the service port we assume that there is a switchover period.
\end{itemize}

In this paper we are interested in the revenue of the system. Every served customer generates  a profit and every lost customer incurs a loss
to the system. Assume that
\begin{itemize}
\item a customer of type $j$ served at station~$i$ gives a profit $\gamma_{ij}$.
\item a customer of type $j$ dropped at station~$i$ causes a penalty $\theta_{ij}$.
\end{itemize}
The motivation for the above assumptions is as follows:
\begin{itemize}
 \item For every packet served the server gains a profit. This profit depends on both the type of packet as well as the source of packet.
 Hence the profit, $\gamma_{ij}$, depends on both $i$ and $j$. 
 \item Further the server has an obligation to meet the contract it has with each source. If the server fails to meet this
 contract it incurs a penalty (loss of packets/reputation/further contracts). This
 again depends on the type of packet and the packet source. Hence the penalty, $\theta_{ij}$, depends on both $i$ and $j$.
\end{itemize}
\section{Performance Measures}
\label{sec:perform}
In this section
we derive expressions for the mean numbers of customers in the stations/retrial loops at various epochs, and use these
to determine an expression for the mean revenue at each station per cycle.
\subsubsection{Mean number of customers at different time epochs}
We know that a communication system works in frame time, where each frame time is fixed. We now assume that the total cycle time
$C$ is this fixed frame time. We have $C= \sum_i(S_i + V_i)$, and $X_{ij}$ and $Y_{ij}$ represent the number of customers
(packets)
of type~$j$ at station~$i$ at the start and end of a visit period of station~$i$ in steady state.

We have
\begin{eqnarray*}
\E[X_{ij}]&=& \E[Y_{ij}] + \lambda_{ij} (C-V_i),\\
\E[Y_{ij}] &=& \E[X_{ij}] (1-r_i(V_i)).\\
\end{eqnarray*}

By solving the above equations we get
\begin{eqnarray*}
\E[X_{ij}]&=&\frac{\lambda_{ij} (C - V_i)}{r_i (V_i)  }, \\
\E[Y_{ij}]&=& \lambda_{ij} (C - V_i)\frac{ 1-r_i (V_i)}{r_i (V_i)}.
\end{eqnarray*}

The customers of type $j$ served during a visit of station~$i$, $T_{ij}$, are the newly arriving customers and the customers in the retrial queues who retry during the visit; hence
\begin{equation}
\E  [T_{ij}]=\E[X_{ij}] p_i(V_i) + \lambda_{ij} V_i.
\label{serve}
 \end{equation}

The customers of type $j$ lost at the end of the visit of station~$i$, $L_{ij}$, are the customers in the retrial queues who did not retry and were dropped at
the end of the visit. Their mean number $\E [L_{ij}]$ is given by:
\begin{equation}
 \E [L_{ij}]=\E[X_{ij}](1- p_i(V_i))q_i(V_i).
\label{lost}
 \end{equation}

\subsubsection{Revenue}

We will now calculate the mean revenue, $R_{i}(V_i)$, after each visit at station~$i$. From Eqs. \eqref{serve}, \eqref{lost}, and the assumption
that a customer of type $j$ served at station~$i$ gives a profit $\gamma_{ij}$ and a customer of type $j$ dropped at station~$i$ causes a penalty $\theta_{ij}$,
we get,
\begin{equation}
 R_{i}(V_i) = \sum_j [\gamma_{ij} \E [T_{ij}] -\theta_{ij} \E [L_{ij}]].
 \label{revenue_i_ini}
\end{equation}

One can alternatively see this as follows. According to our model, all arrivals during the visit time at station~$i$, $V_i$, get served, yielding the profit $V_i \lambda_{ij} \gamma_{ij} $.
Further the arrivals during the non-visit time at station~$i$, $(C-V_i)$, get served with probability $p_i(V_i)/r_i(V_i)$,
which is the conditional probability of a retrial given that the customer disappears during the cycle -- either because of a retrial or because of being dropped.
The server also incurs a loss from the arrivals during the non-visit time at station~$i$, $(C-V_i)$, who are lost with probability $q_i(V_i)/r_i(V_i) $,
which is the conditional probability of being dropped given that the customer disappears during the cycle. Hence
we get,

\footnotesize
\begin{eqnarray}
R_i(V_i) = \sum_j \lambda_{ij} (\gamma_{ij} + \theta_{ij})\left[(C-V_i) \frac{p_{i}(V_i)}{r_i (V_i)}+ V_i\right]- C\sum_j \lambda_{ij}\theta_{ij} .
\label{revenue_i}
\end{eqnarray}
\normalsize
One can verify that Eqs. \eqref{revenue_i_ini}  and \eqref{revenue_i} are the same.

Equation \eqref{revenue_i} can be divided into the $\theta_{ij}$ part and the $\gamma_{ij} + \theta_{ij}$ part.
Hence from the form of the equation we can assume that the system incurs a cost $\theta_{ij}$ for every incoming packet
irrespective of its final state (served or lost) and gains $\gamma_{ij}+ \theta_{ij}$ for every served packet.

\section{Revenue  optimization}
\label{sec:resource}
The system administrator has a limited resource $C$ (frame time) which has to be divided among all the ports.
We assume that the system aims to maximize revenue under the condition of limited available resources, i.e., choose $V_i$
such that $\sum_i R_i(V_i)$
is maximal
while $C=\sum_i (V_i+ S_i)$ is fixed.
 
 Let $\Gamma_i =\sum_j \lambda_{ij} (\gamma_{ij} + \theta_{ij})$ and $M_i(V_i) = R_i(V_i) +C \sum_j \lambda_{ij} \theta_{ij}$. We get
 \begin{equation}
M_i(V_i) = \Gamma_i \left[(C-V_i) \frac{p_{i}(V_i)}{r_i (V_i) }+ V_i\right].
\label{margin_i}
\end{equation}

The maximization of $\sum_i R_i(V_i)$ w.r.t. $V_i$ clearly is the same as the maximization of $\sum_i M_i(V_i)$. Here $M_i(V_i)$ can be interpreted as the
gross profit of the system from station~$i$ and $\Gamma_i$ as the maximum gain per unit time. From here on, we shall also call $M_i(V_i)$ the revenue.

We now have the following optimization problem \\ \textcolor{blue}{\bf REVENUE}
\begin{eqnarray*}
\text{~~~~~~~max}\sum_i M_i(V_i) \\
\text{subject to } \sum_i V_i &=&C -\sum_i S_i\\
\text{and~~~~~~~~~~~~} V_i &\geq& 0 ,~~ \forall i.
\end{eqnarray*}

Differentiating \eqref{margin_i} w.r.t. $V_i$ gives
\footnotesize
\begin{eqnarray}
&M_i ^\prime (V_i)& = \Gamma_i \Bigg[ 1- \frac{p_i(V_i)}{r_i (V_i)} \nonumber \\
&+& (C-V_i)\frac{p^\prime_i(V_i)q_i(V_i)-p_i(V_i)(1-p_i(V_i))q^\prime_i(V_i)}{r_i (V_i) ^2} \Bigg].
\label{margin_diff}
\end{eqnarray}
\normalsize
Further differentiating \eqref{margin_diff} w.r.t. $V_i$ we get

\footnotesize
\begin{eqnarray}
&M_i ^{\prime \prime} (V_i)&= -2\Gamma_i \Bigg[ \frac{p^\prime_i(V_i)q_i(V_i)-p_i(V_i)(1-p_i(V_i))q^\prime_i(V_i)}{r_i (V_i) ^2} \nonumber \\
&+\frac{C-V_i}{ r_i (V_i)^2}&  \Big[-p_i(V_i) p^{\prime}_i(V_i)q^{\prime}_i(V_i)  \nonumber\\
&&+\frac{p_i(V_i)(1-p_i(V_i))q^{\prime \prime}_i(V_i)-p^{\prime \prime}_i(V_i)q_i(V_i)}{2} \nonumber \\
&+ \frac{r_i^{\prime}(V_i)}{r_i (V_i)}&\big[ p^\prime_i(V_i)q_i(V_i)- p_i(V_i)(1-p_i(V_i))q^\prime_i(V_i)\big]\Big]\Bigg].
\label{margin_diff2}
\end{eqnarray}
\normalsize

Here are some logical choices for $p_i(V_i)$ and $q_i(V_i)$.
\begin{itemize}
 \item The longer the visit period the higher the chance a packet will retry. Hence $p_i(V_i)$ can be assumed to be an increasing function in $V_i$. 
 \item The longer the visit period the higher the chance a packet will get served. Hence the packets might be still valuable at the next visit period, which  suggests that
 $q_i(V_i)$ is a decreasing function in $V_i$.
\end{itemize}

Under these assumptions the expression in Eq. \eqref{margin_diff} is readily seen to always be positive which means the revenue obtained from station~$i$ increases with
 the increase in the length of the visit time $V_i$.
In the sequel we shall in addition assume that all $p_i(V_i)$ are concave,
all $q_i(V_i)$ are convex,
and all $r_i(V_i)$ are increasing functions.
These are sufficient conditions for the expression in Eq. \eqref{margin_diff2} to be negative,
so for the objective function to be concave in each of its $N$ components.
Remember that we defined $r_i(V_i)$ as the probability that a customer in the retrial queue leaves the system after the visit period. A reasonable choice for
$q_i(V_i)$ is such that the overall traffic in the buffer decreases with increasing $V_i$. Hence $r_i^{\prime}(V_i) \geq0$ is a reasonable assumption in
many practical situations. 
The above assumptions of concavity and convexity are also quite reasonable, in view of the fact that we consider functions
which are converging to $1$ ($p_i(V_i)$) and $0$ ($q_i(V_i)$), respectively.

The resulting form of optimization problem is widely studied in resource allocation. It is a so-called separable concave optimization problem \cite{Ibaraki_Katoh},  Ch. 2;
the $i$-th term of the objective function only involves $V_i$, and no other $V_j$, $j \neq i$,
and each component is concave.
Such a separable concave optimization problem can be solved
using existing algorithms from \cite{Ibaraki_Katoh}, like RANK.
Without the concavity, one could also solve such separable problems, but the optimization procedure then is much more involved.

We now give a simple step-by-step guideline to follow the RANK procedure outlined in Section 2.2 of \cite{Ibaraki_Katoh}.
Assuming that the functions
$M_i(V_i)$ are concave increasing, we get that the functions $M_i^\prime(V_i)$ are decreasing. Let $C_s = C - \sum_i S_i >0$ represent the total
available time to be divided amongst the stations.
 
\begin{itemize}
\item Calculate all $M_i^{\prime}(0)$ and sort them in decreasing order, say $M_1^{\prime}(0) \geq M_2^{\prime}(0) \geq \cdots \geq M_N^{\prime}(0)$.
\item Allocate total available time $C_s$ to the station with highest slope $M_i^{\prime}(0)$ at $V_i =0$, in our case station~$1$.
 \item Compute $M_1^{\prime}(C_s)$.
 
\begin{itemize}
\item If $M_1^{\prime}(C_s)\geq M_2^{\prime}(0) $ then the procedure stops; optimal strategy is $V_1= C_s$ and $V_i = 0$, $i=2,3,\cdots,N$.
\item If $M_1^{\prime}(C_s) < M_2^{\prime}(0) $, then solve for $U_1$ and $U_2$ such that  $M_1^{\prime}(U_1) = M_2^{\prime}(U_2) $ and 
$U_1 + U_2 =C_s$, i.e., total time is divided between the two stations such that if there is any small additional time available it can be given
to either station~$1$ or station~$2$, giving us the same revenue.
\begin{itemize}
\item If $M_1^{\prime}(U_1)\geq M_3^{\prime}(0) $ then the procedure stops; the optimal strategy is $V_1= U_1$, $V_2= U_2$ and $V_i = 0$, $i=3,\cdots,N$.
\item If $M_1^{\prime}(U_1) < M_3^{\prime}(0) $, then solve for $W_1$, $W_2$ and $W_3$ such that  $M_1^{\prime}(W_1) = M_2^{\prime}(W_2)  = M_3^{\prime}(W_3)$
and $W_1 + W_2 + W_3 =C_s$, i.e., total time is divided between the three stations such that if there is any small additional time available it can be given to either station~$1$ or station~$2$
or station~$3$, giving us the same revenue.
 \item And so on.
\end{itemize}
\end{itemize}
\item  As seen above, the procedure may end with an allocation where some $V_i$ are zero; otherwise after $N$ steps it ends when $C_s$ is allocated amongst all stations. 
\end{itemize}
\section{Numerical examples}
\label{sec:numeric}
In this section we will give two examples of optimal choices of visit times for different stations under some specific conditions on $p_i(V_i)$ and $q_i(V_i)$. 
In each example we assume that irrespective of $V_i$ being positive or not, there is a switchover time $S_i$.
The first, very simple, two-station example is included because it gives insight into the structure of the solution;
in this case one not even needs to use the above-mentioned RANK algorithm.

\subsection*{Example~1}

This example is motivated by current state optical fiber delay loops. Usually in a simple routing node, the delay created by each fiber delay loop is of some fixed length, say $d$.
We assume that the probability of retrial changes linearly with the length of the visit period,
further if the length of the visit period is greater than the length of the delay then all the packets retry and are served:

\begin{eqnarray*}
 p_i(V_i)&=&  \begin{cases}
        V_i/d, & 0 \leq V_i \leq d, \\
     1,  & V_i > d.
      \end{cases}
        \end{eqnarray*}
We further assume that all the packets that are not served in a visit are dropped at the end of it. Hence $q_i(V_i)=1$. 
Now we have,
   \begin{eqnarray*} 
      M_i(V_i) & = & \begin{cases}
        \frac{\Gamma_i V_i(C+d-V_i)}{d }, & 0 \leq V_i \leq d, \\
         \Gamma_i C,  & V_i > d.
      \end{cases}
    \end{eqnarray*}

We solve the above optimization problem REVENUE
with this choice of $M_i(\cdot)$ when $N=2$.
We have $C_s = C-S_1-S_2$.
For the above setting we get 7 different cases under 3 different scenarios.
The first 3 cases represent the scarce resource scenario, i.e., $0<C_s<d$, the next three cases represent limited (but not scarce) scenarios, i.e., $d \leq C_s \leq 2d$, and the last case represents
an abundant resource scenario, i.e., $C_s>2d$.

\begin{enumerate}

 \item When $0<C_s<d$ and $\frac{\Gamma_1}{\Gamma_2}\leq \frac{d+2(S_1+S_2)-C}{C+d}$:
\begin{eqnarray*}
  V_1 &=& 0,\\
   V_2 &=& C_s.
\end{eqnarray*}
 \item When $0<C_s<d$ and $\frac{\Gamma_2}{\Gamma_1}\leq \frac{d+2(S_1+S_2)-C}{C+d}$:
\begin{eqnarray*}
  V_1 &=& C_s,\\
   V_2 &=& 0.
\end{eqnarray*}

\item When $0<C_s<d$, $\frac{\Gamma_1}{\Gamma_2}~>~\frac{d+2(S_1+S_2)-C}{C+d}$\\ and~$\frac{\Gamma_2}{\Gamma_1}~>~\frac{d+2(S_1+S_2)-C}{C+d}$:

\footnotesize
\begin{eqnarray*}
  V_1 &=& \frac{1}{2} \frac{\Gamma_1(C_s+d+S_1+S_2) + \Gamma_2(C_s-d-S_1-S_2)}{\Gamma_1+\Gamma_2},\\
  V_2 &=& \frac{1}{2} \frac{\Gamma_1(C_s-d-S_1-S_2) + \Gamma_2(C_s+d+S_1+S_2)}{\Gamma_1+\Gamma_2}.
 \end{eqnarray*}
\normalsize

\item When $d \leq C_s < 2d$, $\frac{\Gamma_1}{\Gamma_2}~<~\frac{3d+2(S_1+S_2)-C}{C-d}$\\ and~$\frac{\Gamma_2}{\Gamma_1}~<~\frac{3d+2(S_1+S_2)-C}{C-d}$:
 
 \footnotesize
 \begin{eqnarray*}
  V_1 &=& \frac{1}{2} \frac{\Gamma_1(C_s+d+S_1+S_2) + \Gamma_2(C_s-d-S_1-S_2)}{\Gamma_1+\Gamma_2},\\
  V_2 &=& \frac{1}{2} \frac{\Gamma_1(C_s-d-S_1-S_2) + \Gamma_2(C_s+d+S_1+S_2)}{\Gamma_1+\Gamma_2}.
 \end{eqnarray*}
 \normalsize
 
\item When $d \leq C_s < 2d$ and $\frac{\Gamma_1}{\Gamma_2} \geq \frac{3d+2(S_1+S_2)-C}{C-d}$:
\begin{eqnarray*}
 V_1&=& d,\\
 V_2&=& C_s-d.
\end{eqnarray*}
\item When $d \leq C_s < 2d$ and $\frac{\Gamma_2}{\Gamma_1} \geq \frac{3d+2(S_1+S_2)-C}{C-d}$:
\begin{eqnarray*}
 V_1&=& C_s-d,\\
 V_2&=& d.
\end{eqnarray*}

\item When $C_s \geq 2d$:
\begin{eqnarray*}
 V_1&\geq& d,\\
 V_2&\geq& d.
\end{eqnarray*}
Any such combination with $V_1 + V_2 = C_s$ gives the same revenue.
\end{enumerate}

\begin{figure}

 \includegraphics[width=0.5\textwidth, height=0.25\textheight]{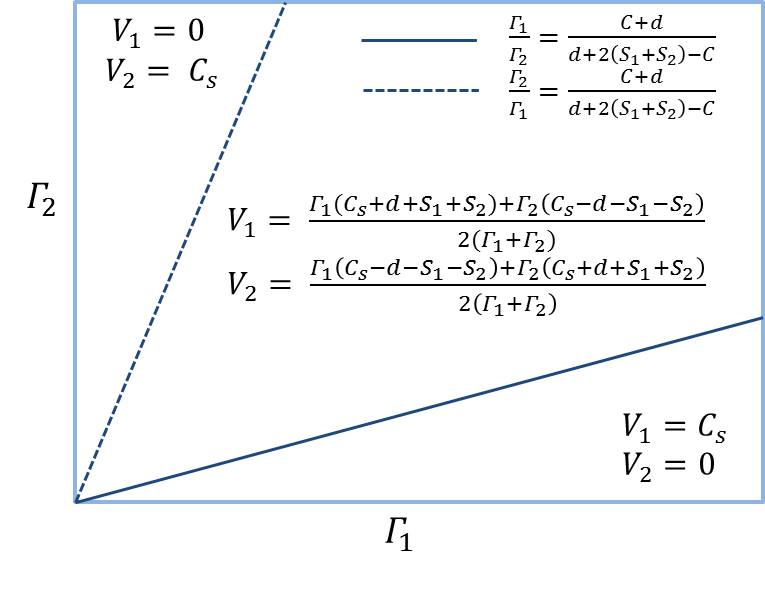} 
 
 \caption{$0<C_s<d$}

 \label{fig_1}
\end{figure} 

\begin{figure}
 \includegraphics[width=0.5\textwidth, height=0.25\textheight]{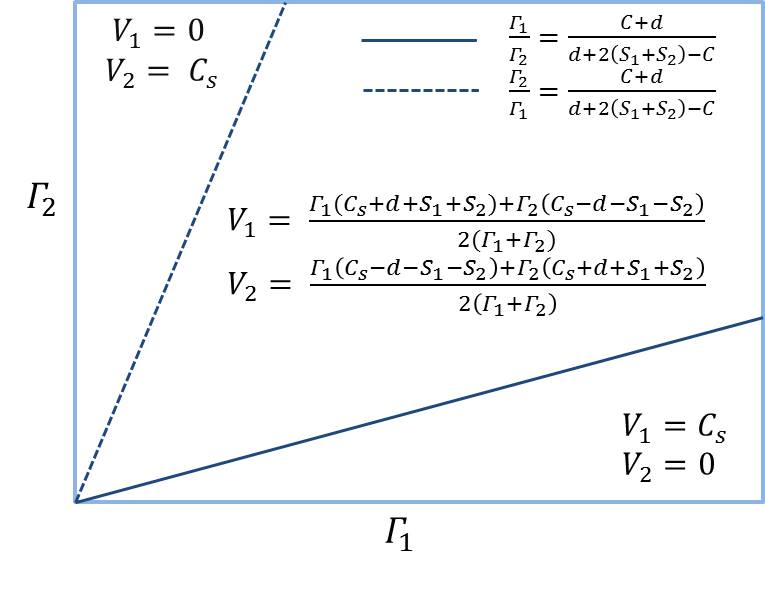} 
 \caption{$d\leq C_s \leq 2d$}
 \label{fig_2}
\end{figure}

For the cases $(1)-(3)$ and the cases $(4)-(6)$, Fig.~\ref{fig_1} and Fig.~\ref{fig_2} respectively
show the optimal choices of
 $V_1$ and $V_2$ in different $\Gamma$ regions.

In this example the $\Gamma_i$ clearly are the key factors which the system administrator should use to make an optimal allocation.

\subsection*{Example 2}
In this example we assume that the packets in the retrial loops of station~$i$ retry after an exponentially distributed time,
with mean $1/\nu_i$. Hence the probability of retrial during $V_i$, $p_i(V_i) = 1 - e^{-\nu_i V_i}$.
Further all the packets which do not retry during $V_i$ are dropped independently with a fixed probability;
$q_i(V_i) = k_i.$	

{\it Case 1:} We consider a 3 station model where all parameters are symmetric except $\Gamma_i$, $i=1,2,3.$ Let $C=14$, $S_i=2$,
$\nu_i=1$, $k_i=0.5$,  $ \forall i=1,2,3.$ In Table \ref{tab:1} we show the optimal values of $V_1,V_2,V_3$ and $\sum_iM_i(V_i)$
for different values of $\Gamma_i$, $i=1,2,3.$ Notice that for $\Gamma_1 = \Gamma_2 = 3$, $V_3$~becomes zero when $\Gamma_3$ drops below $0.01$.

\begin{table}[h]
\resizebox{0.5\textwidth}{!}{
\begin{minipage}{\linewidth}
\begin{center}

\begin{tabular}{| l | c  | c | c | c | c|r|}
     \hline
    $\Gamma_1$ & $\Gamma_2 $ & $\Gamma_3$  & $V_1$ & $V_2$ & $V_3$ & $\sum_iM_i(V_i)$  \\ \hline
    3&3&3&2.6666&2.6666&2.6666&122.3288 \\ \hline    
    3&3&2&2.7837&2.7837&2.4326&108.7920 \\ \hline
    3&3&1&2.9809&2.9809&2.0382&95.4454\\ \hline
    3&3&0.011&3.9949&3.9949&0.0102&83.4455\\ \hline 
    3&3&0.01&4.0000&4.0000&0.0000&83.4455\\ \hline   
    3&2&2&2.9024&2.5483&2.5483&95.1972 \\ \hline
    3&2&1&3.1022&2.7456&2.1522&81.7707\\ \hline
    3&0.01&0.01&5.9308&1.0346&1.0346&42.1918\\ \hline
          \end{tabular} 
          \captionof{table}{Optimal visit length and corresponding maximum revenue for different $\Gamma_i$}
          \label{tab:1} 
\end{center}
      
\end{minipage}}
\end{table}

{\it Case 2:} We consider a 3 station model where all parameters are symmetric except $\nu_i$, $i=1,2,3.$ Let $C=14$, $S_i=2$,
$\Gamma_i=3$, $k_i=0.5$,  $ \forall i=1,2,3.$ In Table \ref{tab:2} we show the optimal values of $V_1,V_2,V_3$ and $\sum_iM_i(V_i)$ for different values of $\nu_i$, $i=1,2,3.$
 
 \begin{table}[h]
\resizebox{0.5\textwidth}{!}{
\begin{minipage}{\linewidth}
\begin{center}

\begin{tabular}{| l | c  | c | c | c | c|r|}
     \hline

     $\nu_1$ & $\nu_2 $ & $\nu_3$  & $V_1$ & $V_2$ & $V_3$ & $\sum_iM_i(V_i)$  \\ \hline

     1&1&1&2.6666&2.6666&2.6666&122.3288 \\ \hline    
    1&1&1.5&2.8959&2.8959&2.2082&123.4510 \\ \hline
    1&1&2&3.0552&3.0552&1.8896&123.9960\\ \hline
    1&1.5&1.5&3.1836&2.4082&2.4082&124.3620\\ \hline 
    0&1.5&1.5&4.8316&1.5842&1.5842&94.8662\\ \hline
          \end{tabular} 
          \captionof{table}{Optimal visit length and corresponding maximum revenue for different $\nu_i$}
          \label{tab:2} 
\end{center}
      
\end{minipage}}
\end{table}

{\it Case 3:} We consider a 3 station model where all parameters are symmetric except $k_i$, $i=1,2,3.$ Let $C=14$, $S_i=2$,
$\nu_i=1$, $\Gamma_i=3$,  $ \forall i=1,2,3.$ In Table \ref{tab:3} we show the optimal values of $V_1,V_2,V_3$ and $\sum_iM_i(V_i)$ for different values of $k_i$, $i=1,2,3.$
 \begin{table}[h]
\resizebox{0.5\textwidth}{!}{
\begin{minipage}{\linewidth}
\begin{center}

\begin{tabular}{| l | c  | c | c | c | c|r|}
     \hline
    $k_1$ & $k_2 $ & $k_3$  & $V_1$ & $V_2$ & $V_3$ & $\sum_iM_i(V_i)$  \\ \hline
    0.5&0.5&0.5&2.6666&2.6666&2.6666&122.3288 \\ \hline    
    0.5&0.5&0.75&2.5568&2.5568&2.8864&121.8160 \\ \hline    
    0.5&0.5&1&2.4784&2.4784&3.0432&121.4070 \\ \hline    
    0.5&1&1&2.3002&2.8500&2.8500&120.2780 \\ \hline    
    0.01&0.5&0.5&0.8730&3.5635&3.5635&124.8190 \\ \hline    
    0.01&1&1&0.6704&3.6648&3.6648&123.9980  \\ \hline    
          \end{tabular} 
          \captionof{table}{Optimal visit length and corresponding maximum revenue for different $k_i$}
          \label{tab:3} 
\end{center}
      
\end{minipage}}
\end{table}

The numerical results suggest that
\begin{itemize}
\item If the value of $\Gamma_i$ increases at a station, then the visit time $V_i$ should increase
(the server serves those stations longer at which it can make a higher profit).
\item If the mean retrial time $1/\nu_i$ decreases at a station, then the visit time $V_i$ should decrease
(the server serves those stations for a shorter time where it can serve many customers in a short span).
\item If the dropping probability $k_i$ increases at a station, then the visit time $V_i$ should increase
(the server serves the stations such that it has fewer lost customers).
\item The optimal $V_i$ not only depend on the parameters of station~$i$ but on parameters of
all the stations.
\end{itemize}

Generally speaking the visit times are chosen such that the system gains higher profit (case~$1$) and provides
better quality (cases $2$ and $3$).

\section{Conclusion}
\label{sec:conclusion}
In this paper we have considered a revenue structure for an optical routing node,
with the aim of providing better Quality-of-Service
to customers of various types.
Modeling an optical routing node as a single server $N$-queue polling system, and demanding that the cycle time of the server is constant,
our goal was to maximize the mean profit per cycle by appropriately choosing
fixed lengths $V_i$ of the visit periods for queue $i$, $i=1,2,\dots,N$.
We have shown that this optimization problem is a separable resource allocation problem which,
under natural assumptions regarding the retrial probability and the dropping probability, becomes
a separable concave resource allocation problem -- a well-studied problem which can be solved using existing algorithms.
We have demonstrated the use of the algorithm RANK for several examples.

In future research we would like to consider the following extensions:
(i) Relax the assumption of cyclic service.
(ii) Relax the assumption of infinite service rate.
(iii) If an optical routing node can use several different wavelengths, then we are faced with a
{\em multi-server} polling system. It will be very interesting to consider the revenue maximization problem
for such a multi-wavelength system.

\section*{Acknowledgment}
The authors gratefully acknowledge a discussion with Cor Hurkens about the algorithm RANK.
The research is supported by the IAP program BESTCOM, funded by the Belgian government,
and by the Gravity program NETWORKS, funded by the Dutch government.

\end{document}